\documentclass[11pt]{article}
\usepackage{graphicx}
\usepackage{amssymb, color}
\usepackage{amsmath}
\def\L{{\bf L}}

\def\D{{\mathcal D}}

\def\ve{\varepsilon}

\def\Ups{\Upsilon}
\def\implies{\Longrightarrow}

\def\Tilde{\widetilde}

\def\C{{\mathcal C}}

\def\bfn{{\bf n}}

\def\bfv{{\bf v}}
\def\bfq{{\bf q}}
\def\bfw{{\bf w}}

\def\bfe{{\bf e}}

\def\curl{\hbox{curl }}

\def\span{\hbox{span}}
\def\ds{\displaystyle}
\def\sqr#1#2{\vbox{\hrule height .#2pt
\hbox{\vrule width .#2pt height #1pt \kern #1pt
\vrule width .#2pt}\hrule height .#2pt }}
\def\square{\sqr74}
\def\endproof{\hphantom{MM}\hfill\llap{$\square$}\goodbreak}
\def\bega{\begin{array}}
\def\enda{\end{array}}
\def\begi{\begin{itemize}}
\def\endi{\end{itemize}}

\def\la{\bigl\langle}
\def\ra{\bigr\rangle}
\def\R{{\mathbb R}}

\def\Supp{\hbox{Supp}}
\def\ov{\overline}
\def\Tilde{\widetilde}
\def\forall{\hbox{for all }~}
\def\v{\vskip 1em}

\def\div{\hbox{div }}
\def\be{\begin{equation}}
\def\beq{\begin{equation}}
\def\bel{\begin{equation}\label}
\def\eeq{\end{equation}}

\newtheorem{theorem}{Theorem}[section]

\newtheorem{lemma}{Lemma}[section]

\newtheorem{remark}{Remark}[section]

\begin{document}
\title{\bf A Posteriori Error Estimates for Self-Similar Solutions to the Euler 
Equations
}
\v

\author{
Alberto Bressan and Wen Shen \\ \, \\
Department of Mathematics, Penn State University.\\
University Park, PA~16802, USA.\\
\\  e-mails:~axb62@psu.edu,  ~wxs27@psu.edu}
\date{Dec 15, 2019}
\maketitle
\begin{abstract} 
The main goal of this paper is to analyze a family of ``simplest possible" initial data for which,
as shown by numerical simulations,
the incompressible Euler equations have multiple solutions.
We take here a first step toward a rigorous validation of these numerical results.
Namely, we consider the system of equations corresponding to a self-similar solution, 
restricted to a bounded domain with smooth boundary. 
Given an approximate solution
obtained via a finite dimensional Galerkin method, 
we establish a posteriori error  bounds on the distance between the numerical 
approximation and the exact solution having the same boundary data.
\end{abstract}
\v
\section{Introduction}
\setcounter{equation}{0}

The flow of a homogeneous, incompressible, non-viscous fluid
in $\R^2$ is modeled by the
Euler equations
\bel{E}
\left\{\begin{array}{rll}
u_t +(u\cdot\nabla) u&=~-\nabla p &\qquad \textrm{(balance of momentum),}\\
\textrm{div}\, u&= ~0 & \qquad\textrm{(incompressibility
condition).}
\end{array}\right.
\end{equation}
Here $u=u(t,x)$ denotes the velocity of the fluid, while the scalar function 
$p$ is a pressure. 
The condition $\div u=0$ implies the existence of a stream 
function $\psi$ such that 
\bel{psi}
u~=~\nabla^\perp\psi,\qquad\qquad (u_1,u_2)~=~(-\psi_{x_2}, \, \psi_{x_1}).\eeq

\begin{figure}[htbp]
\centering
 \includegraphics[scale=0.38]{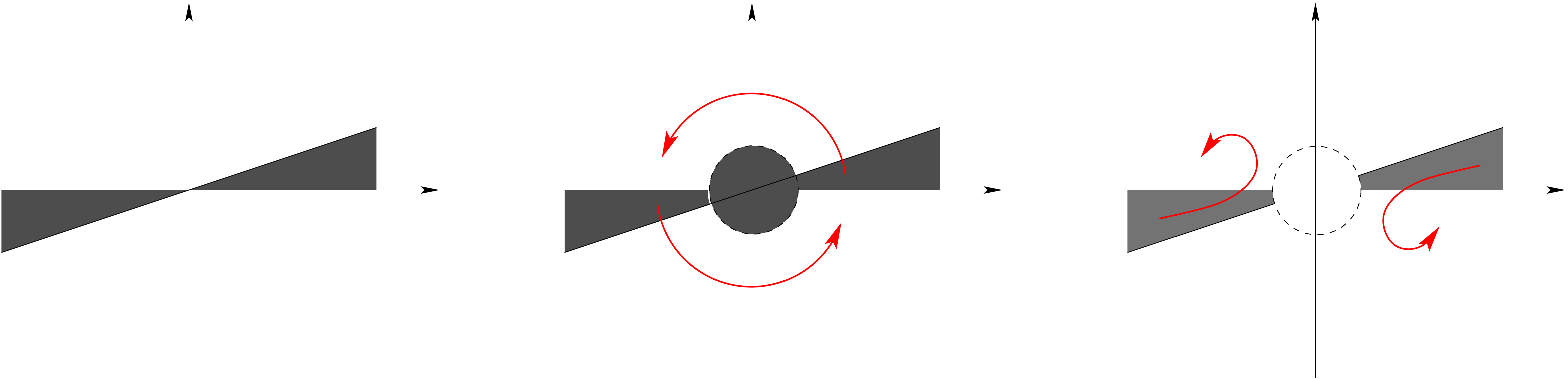}
    \caption{\small  The supports of the initial vorticity considered in (\ref{ssv})-(\ref{ID2}).}
\label{f:e51}
\end{figure}

\begin{figure}[htbp]
\centering
  \includegraphics[scale=0.8]{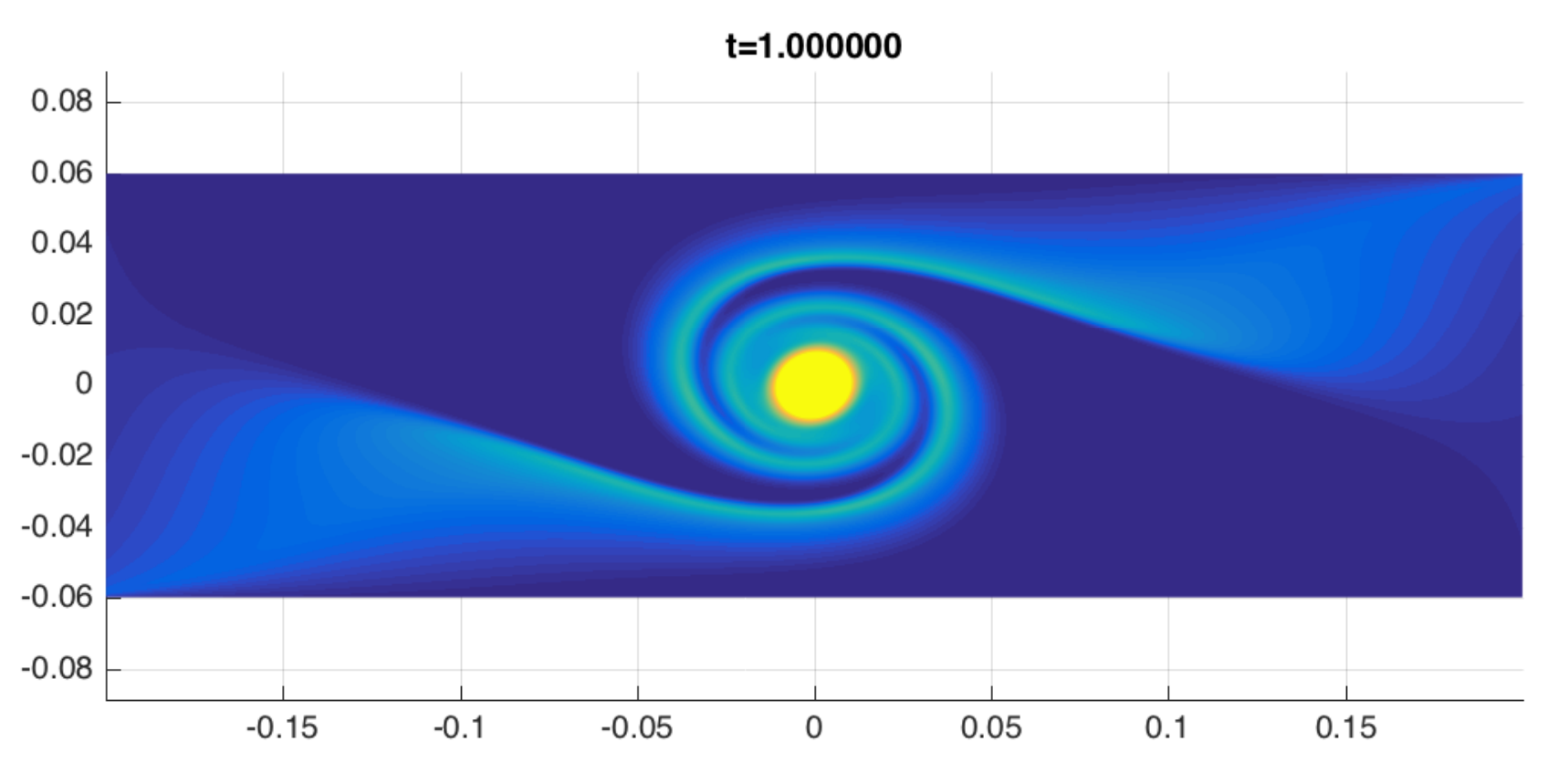}
    \caption{\small 
    The vorticity distribution  at time $t=1$, for a solution to (\ref{E2})
    with initial vorticity $\ov \omega_\ve$.}
\label{f:EC1spi}
\end{figure}

\begin{figure}[htbp]
\centering
  \includegraphics[scale=0.8]{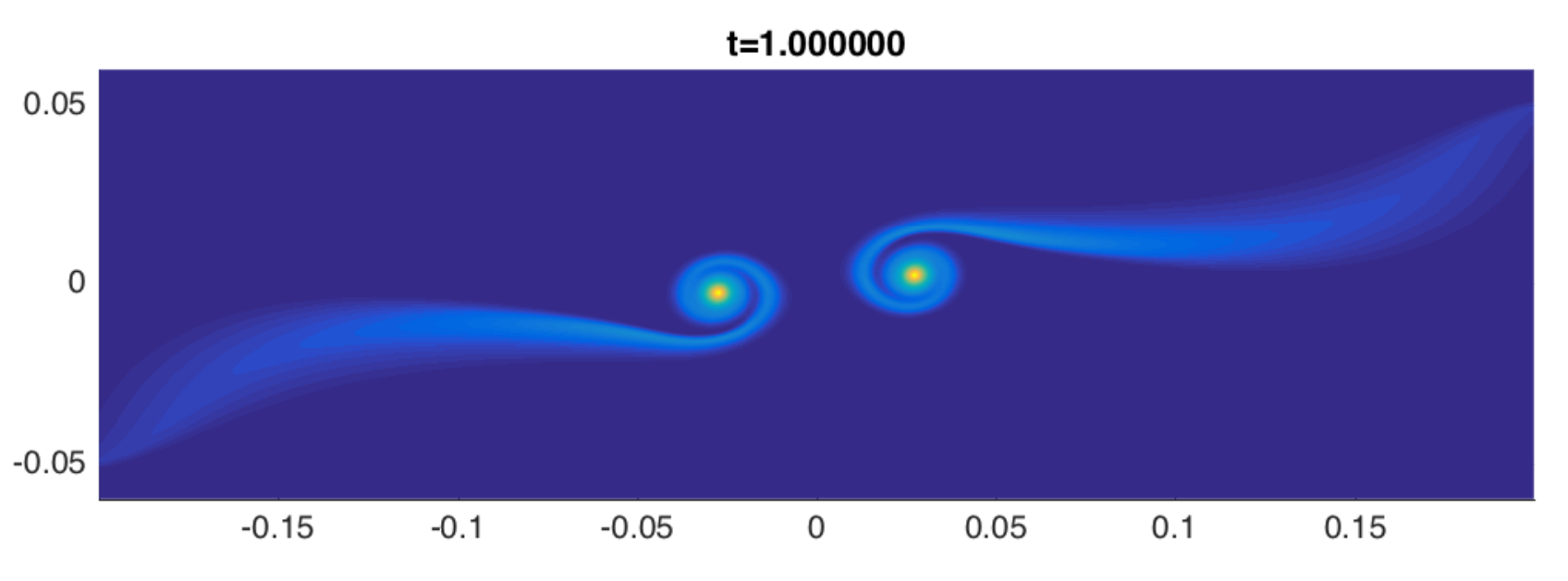}
    \caption{\small
    The vorticity distribution at time $t=1$, for a solution to (\ref{E2})
    with initial vorticity $\ov\omega_\ve^\dagger$.}
\label{f:EC2spi}
\end{figure}

Denoting by $\omega = \curl u = (- u_{1, x_2}+ u_{2, x_1})$ 
the vorticity of the fluid,  
%
it is well known that 
the Euler equations (\ref{E}) can be reformulated 
in terms of the system
\bel{E2}
\left\{\bega{rl} \omega_t + \nabla^\perp \psi\cdot\nabla\omega
&=~0,\\[3mm]
\Delta\psi&=~\omega.\enda\right.\eeq
The velocity $u$ is then recovered from the vorticity
by the Biot-Savart formula
\bel{BS}u(x)~=~{1\over 2\pi}\int_{\R^2} {(x-y)^\perp\over|x-y|^2}
\,\omega(y)\, dy.
\eeq
Our eventual goal is to construct ``simplest possible" initial data
for which the equations (\ref{E}) have multiple solutions. Numerical 
simulations, shown in Figures~\ref{f:EC1spi} and \ref{f:EC2spi}, 
indicate that two distinct solutions can be achieved 
for initial data  where the vorticity 
\bel{om0}\ov \omega(x)~=~ \omega(0,x)~=~\curl  u(0,x)\eeq has the form
\bel{ssv}
\ov \omega(x)~=~r^{-{1/\mu}} \,\ov \Omega(\theta),
\qquad\qquad x= (x_1, x_2) = (r\cos \theta, \,r\sin\theta).\eeq
Here   ${1\over 2}<\mu<+\infty$, while   $\ov\Omega\in \C^\infty(\R)$ is a non-negative, smooth, periodic function
which satisfies
\bel{ovo}
\ov\Omega(\theta) ~=~\ov\Omega(\pi+\theta),\qquad\qquad
\ov\Omega(\theta)~=~0\quad\hbox{if}~~\theta\in \left[{\pi\over 4}\,,\, \pi\right].\eeq
As shown in Figure~\ref{f:e51}, left, the initial vorticity $\ov\omega$ is supported on two wedges, 
and becomes arbitrarily large  as $|x|\to 0$.

One can approximate $\ov\omega$ by two families of initial data
$\ov \omega_\ve, \ov\omega_\ve^\dagger\in \L^\infty(\R^2)$, taking
\bel{ID2}\ov \omega_\ve(x)~=~\left\{\bega{cl} \omega(x)\quad&\hbox{if} ~~|x|>\ve,\cr
\ve^{-1/\mu}\quad&\hbox{if} ~~|x|\leq\ve,\enda\right.\qquad\qquad
\ov \omega^\dagger_\ve(x)~=~\left\{\bega{cl} \omega(x)\quad&\hbox{if} ~~|x|>\ve,\cr
0\quad&\hbox{if} ~~|x|\leq\ve.\enda\right.
\eeq
As $\ve\to 0$, we have $\ov \omega_\ve, \ov\omega_\ve^\dagger \to \ov\omega$
in $\L^p_{loc}$, for a suitable $p$ depending on the parameter  $\mu$ in (\ref{ssv}).
By Yudovich's theorem \cite{Y}, for every $\ve>0$  these initial data yield a unique solution.  However, the numerical simulations indicate that, as $\ve\to 0$,
two distinct limit solutions are obtained. In the first solution, shown in Figure~\ref{f:EC1spi}, both wedges wind up together in a single spiral.   In the second solution, shown in Figure~\ref{f:EC2spi},
each wedge curls up on itself, and two distinct spirals are observed.

\begin{remark} {\rm
The fact that all approximate solutions $\omega_\ve, \omega_\ve^\dagger$ are uniquely determined by their initial data
implies that the ill-posedness exhibited by this example is ``uncurable".   Namely, there is no way to
select one solution with initial datum as in (\ref{ssv})-(\ref{ovo}), preserving the continuous dependence on initial data. }
\end{remark}

\v
We observe that both of these limit solutions are self-similar, 
i.e.~they have the form
\bel{SS3}
\left\{\bega{rl} u(t,x)&=~t^{\mu-1} U\left({x\over t^\mu}\right),\\[4mm]
\omega(t,x)&=~t^{-1} \Omega\left({x\over t^\mu}\right),\\[4mm]
\psi(t,x)&=~t^{2\mu-1} \Psi \left({x\over t^\mu}\right)
 .\enda\right.\eeq
Notice that, by self-similarity, these solutions are completely determined as soon as
we know their values  at time $t=1$.  Indeed, these are given by $U, \Omega,\Psi$.
Inserting (\ref{SS3}) in (\ref{E2}), one obtains the equations
\bel{SSE}\left\{
\bega{rl}\Big(\nabla^\perp \Psi -\mu y\Big)\cdot \nabla\Omega &=~ \Omega\,,\\[3mm]
\Delta \Psi&=~\Omega\,,
\enda
\right.
\eeq
while the velocity is recovered by
\bel{U}U~=~\nabla^\perp \Psi.\eeq

Constructing two distinct self-similar solutions of (\ref{E2}) with the same initial data (\ref{ssv})
amounts to finding two distinct solutions $(\Omega,\Psi)$, $(\Omega^\dagger,
\Psi^\dagger) $ of (\ref{SSE}) with the same 
asymptotic behavior as $|x|\to + \infty$. More precisely, writing the vorticity $\Omega$ in polar coordinates, this 
means
\bel{asy1}\lim_{r\to +\infty} r^{1\over\mu}\,\Omega(r,\theta)~=~\lim_{r\to +\infty} r^{1\over\mu}\,\Omega^\dagger(r,\theta)
~\doteq~\ov\Omega(\theta),\eeq 
for some smooth function $\ov\Omega$ as in (\ref{ovo}).

\begin{figure}[htbp]
\centering
  \includegraphics[scale=0.22]{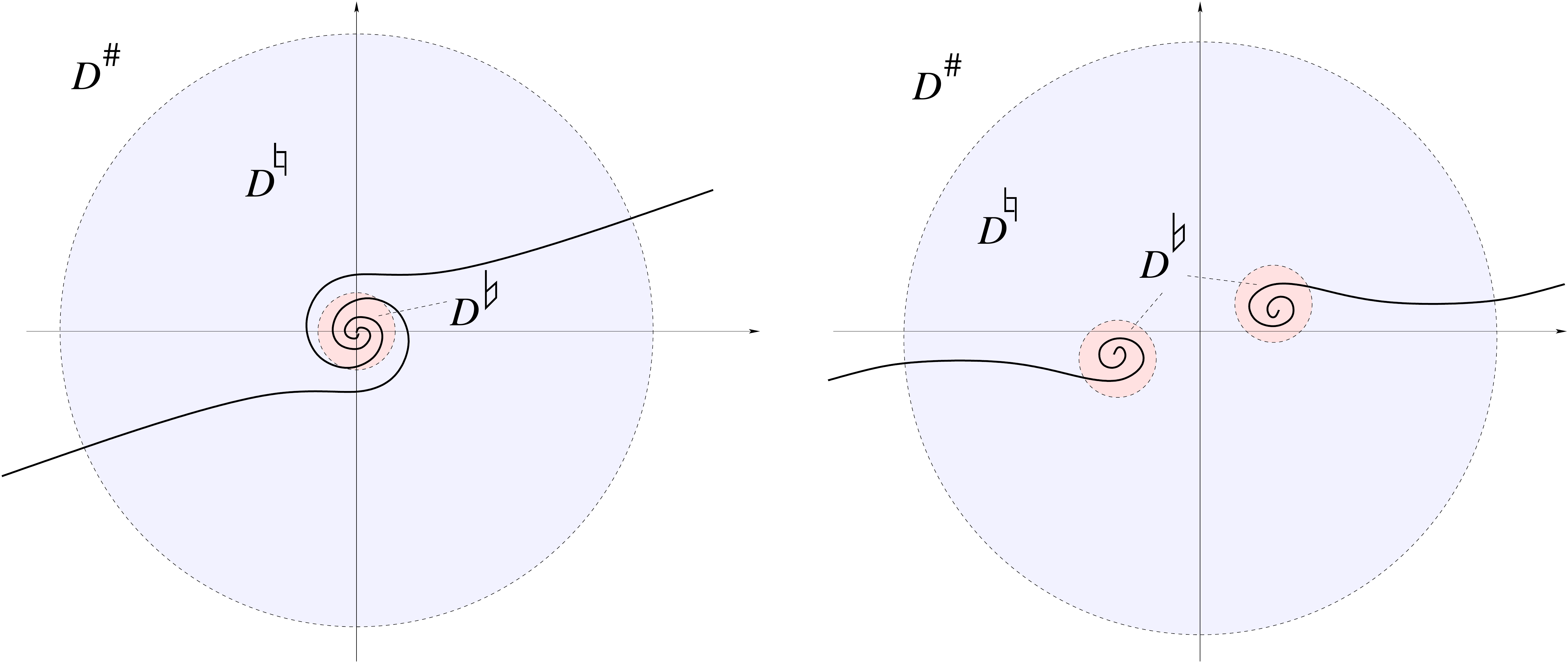}
    \caption{\small Decomposing the plane $\R^2=\D^\sharp\cup\D^\natural\cup\D^\flat$
    into an outer, a middle, and an inner domain. Left: the case of a single spiraling vortex, as in Fig.~\ref{f:EC1spi}. Right: the case of two spiraling vortices, as in Fig.~\ref{f:EC2spi}.}
\label{f:e64}
\end{figure}

Since the two solutions in Figures~\ref{f:EC1spi} and \ref{f:EC2spi}
 are produced by numerical computations, 
a natural question is whether an exact self-similar solution of the 
Euler equations exists, close to each computed one.
This requires suitable a posteriori error bounds.
Toward this goal, two difficulties arise:
\begi
\item[(i)] The self similar solution $(\Omega,\Psi)$ is defined on the entire plane
$\R^2$, while a numerical solution is computed only on some bounded domain.
\item[(ii)] The solution is  smooth,  with the exception of one or two
points corresponding to the spirals' centers.   In a neighborhood of these points
the standard error estimates  break down.
\endi
To address these issues, we propose a domain decomposition method.
As shown in Fig.~\ref{f:e64},
the  plane can be decomposed 
into an outer domain $D^\sharp\supseteq\{x\in\R^2\,; |x|>R\}$, 
an inner domain $\D^\flat$ containing
a neighborhood of the spirals' centers where the solution has singularities, and a bounded intermediate domain $\D^\natural$ where the solution is smooth.
The solution is constructed analytically on $\D^\sharp$ and on $\D^\flat$,
and numerically on $\D^\natural$.   
These three  components are then patched together 
by suitable matching conditions.   

A detailed analysis of the solution to (\ref{SSE}) in a neighborhood of 
infinity and near the spirals' centers will appear in the companion paper
\cite{BM}, relying on the approach developed in  \cite{E1, E2, E3}. 
In the present paper we focus on the derivation of a posteriori
error estimates for a numerically computed solution to (\ref{SSE}), on 
a bounded domain $\D\subset\R^2$ with smooth boundary $\partial\D$.    
As shown in Fig.~\ref{f:e88}, we assume that this boundary
can be decomposed as the union of two closed, disjoint components:
\bel{dbd}\partial\D~=~\Sigma_1\cup\Sigma_2\,.\eeq
We seek a solution to (\ref{SSE}), satisfying boundary conditions of the form
\bel{BC}\left\{
\bega{rl}\Psi(x)&=~g(x),\qquad\qquad x\in \partial \D,\\[3mm]
\Omega(x)&=~h(x),\qquad\qquad x\in\Sigma_1\,,\enda\right.\eeq
where $g,h$ are given smooth functions.
Given an approximate solution 
$(\Omega_0,\Psi_0)\in \C^{0,\alpha}(\D)\times \C^{2,\alpha}(\D)$, computed
by a Galerkin finite dimensional approximation,  
we want to prove the existence
of an exact solution close to the approximate one.

The remainder of the paper is organized as follows. Section~\ref{s:2} introduces the basic framework,
specifying the main assumptions on the numerical scheme and on the approximate solution.
Section~\ref{s:3} begins by  analyzing the first equation in (\ref{lam1}), regarded as a linear 
PDE for the vorticity function $\Omega$.  In this direction, Lemma~\ref{l:42}
provides a detailed estimate on how the solution depends on the vector field $\nabla^\perp \Phi$.
In addition, Lemma~\ref{l:interp} 
yields a sharper regularity estimate on the solutions, deriving an
a priori bound on their $\C^{0,\alpha}$ norm.

Finally, in Section~\ref{s:4} the exact solution $\Omega$ is constructed as the fixed point of a transformation
which is contractive w.r.t.~a norm equivalent to the $\L^2$ norm. We remark that,
in order  to achieve this contractivity, a bound on the norms $\|\Omega\|_{\L^2}$ and $\|\Phi\|_{H^2}$ 
is not good enough. Indeed, we need an a priori bound on $\|\Omega\|_{\C^{0,\alpha}}$ and on 
$\|\Phi\|_{\C^2}$. 
This is achieved by means of Lemma~\ref{l:interp}.
At the end of Section~\ref{s:4} we collect all the various constants appearing in the estimates,
and summarize our analysis  by stating a theorem on the existence of an exact solution, close to the 
computed one.

For results on the uniqueness of solutions to the incompressible Euler equations we refer to 
\cite{BH, MP, Y}.   Examples showing the non-uniqueness of solutions to the incompressible Euler equations 
were first constructed   in \cite{Sch, Shn}.  See also 
\cite{V} for a different approach, yielding more regular solutions. 
Following the major breakthrough in 
\cite{DS09}  several examples 
of multiple solutions for Euler's equations have recently been provided \cite{Dan, DRS, DanSz, DS10}. 
These solutions are obtained by means of convex integration and a Baire category argument.
They have turbulent nature and their physical significance is unclear.  

Our numerical simulations, on the other hand, suggest that ``uncurable"  non-uniqueness
can arise from quite simple initial data. In our case, both solutions can be easily visualized;
they remain everywhere smooth (with the exception of one or two points), and conserve energy at all times.

\section{Setting of the problem}
\setcounter{equation}{0}
\label{s:2}
Throughout the following, we consider the boundary value problem (\ref{SSE}), (\ref{BC}) on 
a bounded, open domain $\D\subset \R^2$, with smooth boundary $\partial\D$ decomposed as in (\ref{dbd}).   
The unit outer normal 
at the boundary point $y\in \partial \D$ will be denoted by $\bfn(y)$.
We call $\Psi^g$ the solution to the  non-homogeneous boundary value problem
\bel{Pg}\left\{\bega{rll}
\Delta\Psi(x)&=~0\qquad & x\in\D,\\[3mm]
\Psi(x)&=~g(x)\qquad &x\in\partial\D,\enda\right.\eeq
and define the vector field
\bel{bbv}
\bfv(x)~\doteq~\nabla^\perp \Psi^g(x) -\mu x\,.\eeq
Assuming that $g$ is smooth, the same is true of $\bfv$.
Given $\Omega\in \L^2(\D)$, we define 
\bel{D-1}\Phi~=~\Delta^{-1}\Omega\eeq 
to be the solution of 
\bel{lam2}
\left\{\bega{rll} \Delta\Phi (x)&=~\Omega(x),\qquad &x\in \D,\\[3mm]
\Phi(x)&=~0,\qquad &x\in \partial \D.\enda\right.\eeq
Our problem thus amounts to finding a function $\Omega$ such that
\bel{lam1}
\left\{\bega{rll} \bigl(\nabla^\perp \Phi(x)+\bfv(x)\bigr)\cdot \nabla\Omega(x)&=~\Omega(x),\qquad\qquad 
&x\in \D,\\[3mm]
\Omega(x)&=~h(x), &x\in \Sigma_1\,,\enda\right.\eeq
where $\Phi=\Delta^{-1}\Omega$ and $\bfv$ is the vector field at (\ref{bbv}).
In the following, for any given function $\Phi\in \C^{2}(\D)$, 
we denote by $\Omega=\Gamma(\Phi)$ the solution to (\ref{lam1}).
We seek a fixed point of the composed map
\bel{ldef}
\Omega~\mapsto~\Lambda(\Omega)~\doteq~\Gamma(\Delta^{-1}\Omega).\eeq
\v
We shall consider approximate solutions of (\ref{lam1}) which are obtained
by finite dimensional Galerkin
approximations.   More precisely, given a finite set of linearly independent functions
$\{\phi_j\,;~~1\leq j\leq N\}\subset\L^2(\D)$, consider  the orthogonal decomposition
$\L^2(\D)~=~U\times V$, where
\bel{UV} U~=~\span\{\phi_1,\ldots, \phi_N\},\qquad \qquad V= U^\perp,\eeq
with orthogonal projections
\bel{proj}
P:\L^2(\D)~\mapsto~ U,\qquad\quad (I-P): \L^2(\D)~\mapsto ~V.\eeq
\v
{\bf Example 1.} We can choose the functions
$\phi_j$ to be piecewise affine, obtained from a
triangulation of the domain $\D$. This implies
\bel{pjp}
\phi_j\in W^{1,\infty}(\D)\subset \C^{0,\alpha}(\D)\eeq
for every $0<\alpha\leq 1$.   In this case, the gradient $\nabla\phi_j\in\L^\infty\cap BV$ 
is piecewise constant, with jumps
along finitely many segments.
\v
{\bf Example 2.} In alternative, one can choose $\phi_1,\ldots,\phi_N\in \L^2(\D)$
to be the  first $N$ normalized eigenfunctions of the Laplacian.   More precisely, 
for $j=1,\ldots,N$ we require that the functions $\phi_j$ satisfy
\bel{EL}\left\{\bega{rl} \Delta\phi_j + \lambda_j \phi_j~=~0\qquad & x\in \D,\\[3mm]
\phi_j~=~0\qquad & x\in \Sigma_1\,,\\[3mm]
\bfn\cdot \nabla\phi_j~=~0\qquad & x\in \Sigma_2\,,\enda\right.\eeq
with eigenvalues $0<\lambda_1\leq \lambda_2\leq\cdots\leq \lambda_N$. Moreover, $\|\phi_j\|_{\L^2}=1$.
\v
\begin{remark}\label{r:1}{\rm In both of the above cases, there may be no linear combination $\sum_j c_j \phi_j$ 
of the basis functions that matches  the boundary data $h$ along $\Sigma_1$. 
This issue can be addressed simply by adding to our basis an additional function $\phi_0\in \C^\infty(\D)$,
chosen so that $\phi_0=h$ on $\Sigma_1$.}
\end{remark}

Our basic question can be formulated as follows. 

\begi
\item[{\bf (Q)}]
{\it 
Assume we can find a finite dimensional approximation
\bel{UU} \Omega_0 ~=~\sum_{j=1}^N c_j\phi_j\,,\qquad\qquad
\Phi_0~=~\Delta^{-1}\Omega_0\,,\eeq
with error
\bel{AA} \Big\| \Omega_0 -P\,\Lambda(\Omega_0)\Big\|_{\L^2(\D)}~=~\delta_0\,.\eeq
How small should $\delta_0$ be, to make sure that an exact solution 
$(\Omega,\Phi)$ of (\ref{lam1}), (\ref{D-1}) exists,
close to  $(\Omega_0,\Phi_0)$ ?}
\endi
\v

Given a function $\Phi\in \C^2(\D)$, the linear, first order PDE  (\ref{lam1}) for $\Omega$ can be solved 
by the method of characteristics.
Namely, consider the
 vector field
\bel{qd}\bfq(x)~=~\nabla^\perp \Phi(x)+\bfv(x)
,\eeq
 whose divergence is
 \bel{divq}
 \div \bfq~=~-2\mu\,.\eeq
We shall denote by $t\mapsto \exp(t\bfq)(y)$ the solution to the ODE
$$\dot x(t)~=~\bfq(x(t)),\qquad x(0)=y.$$
For convenience, we shall use the notation
\bel{xty}t~\mapsto ~x(t,y)~\doteq~\exp(-t\bfq)(y)\eeq
for the solution to the ODE
$$\dot x(t)~=~-\bfq(x(t)),\qquad x(0)=y.$$
Consider the set
\bel{D*} \D^*~\doteq~\Big\{ y\in \ov\D\,;~~x(\tau,y)\in \Sigma_1\cap\Supp(h)\quad\hbox{for some}~\tau\geq 0\Big\}.\eeq
In other words, $y\in \D^*$ if the characteristic through $y$ reaches a boundary point in the support of $h$ within 
finite time.  Calling 
\bel{tauy}\tau(y)~\doteq~\min\,\bigl\{t\geq 0\,;~~ x(t,y)\in \Sigma_1\bigr\}\eeq
the first time when the characteristic starting at $y$ reaches the boundary $\Sigma_1$,
the solution to (\ref{lam1}) is computed by
\bel{OR}\Omega(y)~=~\left\{ \bega{cl} e^{\tau(y)}\, h(x(\tau(y), y))\qquad &\hbox{if} ~~x\in \D^*,\\[3mm]
0\qquad&\hbox{if}~~x\notin\D^*.\enda\right.\eeq

\begin{figure}[htbp]
\centering
  \includegraphics[scale=0.45]{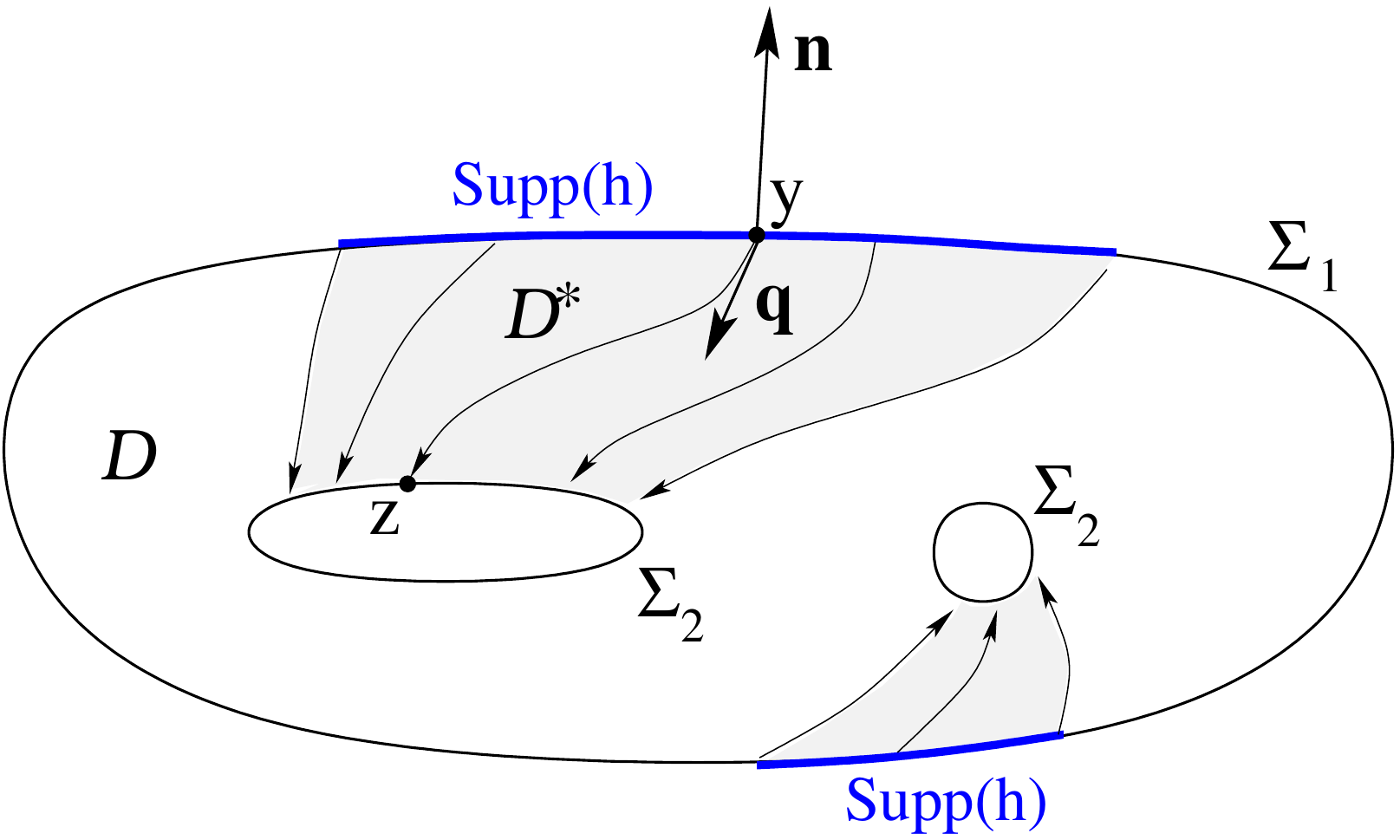}
    \caption{\small  According to {\bf (A1)}, every characteristic starting at a point 
    $y\in \Sigma_1\cap {\rm Supp }(h)$ exits from the domain $\D$ at some boundary point 
    $z=z(y)\in \Sigma_2$, at a time $T(y)\leq T^*$. The shaded region represents the subdomain $\D^*$ in
    (\ref{D*}).}
\label{f:e88}
\end{figure}

The following transversality assumption will play a key role in the sequel (see Fig.~\ref{f:e88}).
\v
\begi
\item[{\bf (A1)}]  {\it There exists  constants $T^*, c_1>0$ such that, for 
every boundary point $y\in \Sigma_1\cap \Supp(h)$,   the following holds.
\begi
\item[(i)] The vector  $\bfq$ 
is strictly inward pointing:
\bel{OT3}\bigl\langle \bfn(y)\,,\, \bfq(y)\bigr\rangle~\leq~-c_1\,.
\eeq
\item[(ii)] The characteristic $t\mapsto \exp(t\bfq)(y)$ remains inside $\D$ 
 until it reaches a boundary point
 \bel{exitp}z(y)\,=\,  \exp\bigl(T(y)\bfq\bigr)(y)~\in ~\Sigma_2\eeq
 within a finite time $T(y)\leq T^*$, and exits transversally:
 \bel{OT2}
\Big\langle \bfn(z(y))\,,\, \bfq(z(y))\Big\rangle~\geq ~c_1\,.\eeq
\endi
}
\endi

As in (\ref{UV})-(\ref{proj}), we consider the decomposition $H\doteq \L^2(\D)=U\times V$, with perpendicular projections 
$P$ and $ I-P$, and write $\Omega= (u,v)$. The partial derivatives of the map
$\Lambda=\Lambda(u,v)$ introduced at (\ref{ldef}) will be denoted by 
$D_u\Lambda$, $D_v\Lambda$.
Let a finite dimensional approximate solution $\Omega_0= (u_0,0)\in U$ of (\ref{lam1}) be given,
with $\Phi_0=\Delta^{-1}\Omega_0$.    Throughout the following, we denote by
\bel{Adef} A~\doteq~P\circ D_u\Lambda(\Omega_0).\eeq
 the  partial differential w.r.t.~$u$ of the map $(u,v)\mapsto P\Lambda(u,v)$, computed at  the point
 $\Omega_0$.
Notice that $A$ is a linear map from the finite dimensional space $U$ into itself.

 Since we are seeking a fixed point, in the same 
spirit of the Newton-Kantorovich theorem~\cite{Berger, Ciarlet, Deimling} 
we shall assume the invertibility of the 
Jacobian map, restricted to the finite dimensional subspace $U$.
\begi
\item[{\bf (A2)}] {\it The operator $I - A$
from $U$ into itself has a bounded inverse, with }
\bel{IA}\big\|(I-A)^{-1}\bigr\|_{{\cal L}(U)}~\leq~\gamma~<~+\infty\eeq
for some $\gamma\geq 1$.
\endi
Here the left hand side refers to the operator norm, in the space of linear operators on $U\subset\L^2(\D)$.
Notice that (\ref{IA})  implies that the operator $I-AP:H\mapsto H$ is also invertible, 
with 
\bel{IA2}
\big\|(I-A\,P)^{-1}\bigr\|_{{\cal L}(H)}~\leq~\gamma.\eeq


Concerning the $V$-component, a  key assumption used in our analysis will be
\begi
\item[{\bf (A3)}] {\it The orthogonal spaces $U,V$ in (\ref{UV}) are chosen so that 
\bel{LVsmall}
\| \Delta^{-1}\circ (I-P)\|~<~\ve_0~<\!\!<~1\,.\eeq
}
\endi
Intuitively this means that, in the decomposition $\Omega = (u,v)\in U\times V$, the 
component $v\in V$ captures the high frequency modes, 
which are heavily  damped by the inverse Laplace operator \cite{BRS, Cle, KS}.
\v
{}From an abstract point of view,  proving the existence of a fixed point of the map $\Omega\mapsto 
\Lambda(\Omega)$ in (\ref{ldef}) is a simple matter.
Using the decomposition $\Omega = (u,v)\in U\times V$ with orthogonal projections as in (\ref{proj}),
we start with an initial guess $\Omega_0 = (u_0,0)$.
Assuming  (\ref{IA}), we write (\ref{lam1}) in the form
$$\Omega-AP\,\Omega~=~\Lambda(\Omega)- AP\,\Omega\,.$$
Equivalently,
\bel{eqeq}
\Omega~=~\Ups(\Omega)~\doteq~(I-AP)^{-1}\bigl(\Lambda(\Omega)- AP\,\Omega\bigr).\eeq
The heart of the matter is to show that, on a neighborhood ${\cal N}$ of $\Omega_0$,
 the map $\Omega\mapsto\Ups(\Omega)$ is a strict contraction. 
 If the initial error $\|\Ups(\Omega_0)-\Omega_0\|$
 is sufficiently small, the iterates $\Omega_n=\Ups^n(\Omega_0)$ will thus remain inside ${\cal N}$ 
 and converge to a fixed point.  The contraction property is proved as follows.
By (\ref{IA2}) one has $\|(I-AP)^{-1}\|\leq\gamma$.   On the other hand, computing the differential of the map
$\Omega\mapsto \Lambda(\Omega)- AP\,\Omega$
w.r.t.~the components $(u,v)\in U\times V$,  at the point $\Omega=\Omega_0= (u_0,0)$  one obtains
\bel{diff}
D\bigl(\Lambda(u,v)-Au\bigr)
~=~\left(
\bega{ccc} 0 && P D_v\Lambda\\[4mm]
(I-P) D_u\Lambda && (I-P) D_v\Lambda\enda\right).\eeq
Because of (\ref{LVsmall}), we expect 
\bel{dism}\|D_v\Lambda\|~\leq~\|D\Gamma\|\cdot \bigl\| \Delta^{-1} \circ(I-P)\bigr\|~<\!<~1.\eeq
By possibly using an equivalent norm on the product space $U\times V$ (see (\ref{n*}) for details), 
we thus achieve the strict contractivity of the map $\Ups$ in (\ref{eqeq}).  
\v
In the next section, more careful estimates will be derived on the differentials of $\Gamma$ and 
$\Lambda$  in (\ref{ldef}). In this direction we remark that,
while the operator $\Delta^{-1}$ is well defined on $\L^2(\D)$, 
a difficulty arises in connection with the differential of $\Gamma$.
Indeed, to compute this differential, we need to perturb the function $\Phi$ in (\ref{lam1})
and estimate the change in the corresponding solution $\Omega$.
This can be done assuming that  $\Phi\in\C^2$, 
hence $\nabla^\perp\Phi\in \C^1$.  Unfortunately, the assumption $\Omega\in \L^2$ 
only implies $\Phi= \Delta^{-1}\Omega\in H^2$, which does not yield any bound on $\|\Phi\|_{\C^2}$.
In order to establish the desired a posteriori error bound,
an additional argument will thus be needed, showing that our approximate solutions actually
enjoy some additional regularity.

\section{Preliminary lemmas}
\setcounter{equation}{0}
\label{s:3}


\begin{figure}[htbp]
\centering
  \includegraphics[scale=0.5]{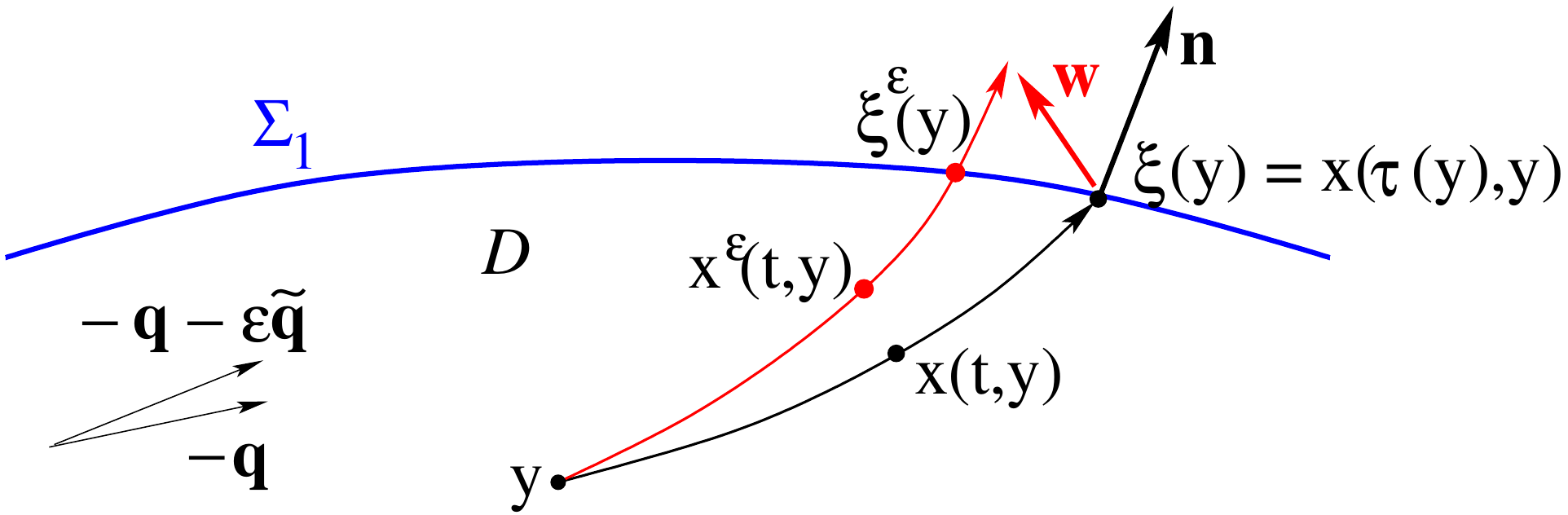}
    \caption{\small Computing the perturbed solution $\Omega^\ve(y)$ in (\ref{om4}), by estimating the change in  
    the characteristic through $y$.}
\label{f:e86}
\end{figure}

\subsection{Continuous dependence of solutions to a first order linear PDE.}
%
As remarked in Section~\ref{s:2}, solutions to the linear PDE (\ref{lam1}) can be found by the method of characteristics (\ref{OR}).  Our present goal is to understand how the solution changes,
depending on the vector field $\bfq$ in (\ref{qd}).  To fix the ideas,
consider the boundary value problem
\bel{lam3}
\left\{\bega{rll} \bfq(x)\cdot \nabla\Omega&=~\Omega,\qquad\qquad 
&x\in \D,\\[3mm]
\Omega&=~h, &x\in \Sigma_1\,,\enda\right.\eeq
assuming that $\bfq:\D\mapsto\R^2$ is a $\C^1$ vector field
satisfying {\bf (A1)} together with 
\bel{qass}\|\bfq\|_{\C^1(\D)}~\leq~M,\qquad\qquad \div \bfq(x) ~=~ -2\mu\qquad\forall x\in \D\,.\eeq
Given a second vector field $\Tilde\bfq\in \C^1$, consider the family of perturbations
\bel{qep}\bfq^\ve(x)~=~\bfq(x) + \ve \,\Tilde\bfq (x),\eeq
and  let 
\bel{om4}\Omega^\ve(x)~=~\Omega(x) +\ve\, \Tilde\Omega (x) + o(\ve)\eeq
be the corresponding solutions of (\ref{lam3}).
Here and in the sequel, the notation $o(\ve)$ indicates a higher order infinitesimal, 
so that $\ve^{-1}o(\ve)\to 0$ as $\ve\to 0$.
The next lemma provides an $\L^2$~estimate on size of the first order perturbation $\Tilde\Omega$.
Setting $\Sigma^*_1\doteq\Sigma_1\cap\Supp(h)$, we introduce the constant
\bel{TC}
\Tilde C~=~ \sup_{x\in \Sigma_1^*} {1\over | \langle \bfn(x), \bfq(x)\rangle|} 
 \cdot \|h\|_{\C^0}+\left( 1+\,\sup_{x\in \Sigma_1^*} {|\bfq(x)|\over | \langle \bfn(x), \bfq(x)\rangle|}\right) \|Dh\|_{\C^0} \,,
 \eeq
and define
\bel{pera} K(M, t)~=~e^t\left({e^{(2M+1)t}- e^{-2\mu t}\over 4M+4\mu+2}\right)^{1/2}
.\eeq

\begin{lemma}\label{l:42} Recalling (\ref{D*}), assume that $\bfq\in \C^1(\D)$ satisfies  (\ref{qass})
together with {\bf (A1)}.
In the  setting considered at (\ref{qep})-(\ref{om4}), the first order perturbation $\Tilde \Omega$ satisfies
\bel{oper}
\|\Tilde\Omega \|_{\L^2(\D)}~\leq~
\Tilde C \cdot K\Big(\|D\bfq\|_{\C^0(\D^*)},\,T^*\Big)\cdot \|\Tilde\bfq\|_{\L^2(\D^*)}\,.\eeq
\end{lemma}

{\bf Proof.} {\bf 1.}  In analogy with (\ref{xty}),
we denote by  
$$ t\,\mapsto\, x^\ve(t,y)\,=\,\exp(-t\bfq^\ve)(y)$$
the solution to
\bel{od2}\dot x^\ve~=~-\bfq(x)-\ve \Tilde\bfq(x),\qquad x(0)=y.\eeq
For $0\leq t<\tau(y)$, the tangent vector
\bel{tanv}\bfw(t,y)~\doteq~\lim_{\ve\to 0} {x^\ve(t,y)-x(t,y)\over\ve}\eeq
provides a solution to the linearized equation
\bel{leq}
\dot \bfw(t,y)~=~-D\bfq(x(t,y))\cdot \bfw(t,y)- \Tilde\bfq(x(t,y)),\qquad\qquad \bfw(0,y)=0.\eeq
For notational convenience, 
we extend the definition of the vector $\bfw(t,y)$ by setting
\bel{tT}\bfw(t,y)~=~\bfw(\tau(y), y)\qquad\hbox{if}\qquad t\in [\tau(y), \,T^*].\eeq
We denote by
\bel{xie}\bega{l} \xi(y)\,=\,x(\tau(y),y)\,=\,\exp\bigl(-\tau(y)\bfq\bigr)(y),
\\[4mm]  \xi^\ve(y)\,=\,x^\ve(\tau^\ve(y),y)\, =\,\exp\bigl(-\tau^\ve(y)\bfq^\ve\bigr)(y)\,,\enda\eeq
the points where the characteristic through $y$ crosses the boundary $\Sigma_1$,
and consider the expansions
\bel{txx}
\xi^\ve(y)~=~\xi(y)+\ve\tilde\xi(y)+o(\ve),\qquad\qquad \tau^\ve(y)~=~\tau(y)+\ve \Tilde\tau(y)+o(\ve).\eeq
Observing that 
\bel{wwy}
\Big\langle \bfn\bigl(x(\tau(y),y)\bigr)~,~\bfw(\tau(y),y) -\Tilde\tau(y)\cdot \bfq\bigl(x(\tau(y),y)\bigr)
\Big\rangle~=~0,
\eeq
we obtain 
\bel{ttau}\Tilde\tau(y)~=~{\Big\langle \bfn\bigl(\xi(y)\bigr)~,~\bfw(\tau(y),y)\Big\rangle \over 
\Big\langle \bfn\bigl(\xi(y)\bigr)~,~\bfq\bigl(\xi(y)\bigr)\Big\rangle}\,,\eeq
\bel{txi}
\tilde \xi(y)~=~\bfw(\tau(y),y) -\Tilde\tau(y)\, \bfq\bigl(\xi(y)\bigr)
\,,\eeq
\bel{tx2}
|\tilde \xi(y)|~\leq~|\bfw(\tau(y),y)|\cdot \left( 1+  \left| { \bfq(\xi(y))\over  \la \bfn(\xi(y)),\,\bfq(\xi(y))
\ra}\right| \,\right). \eeq

Finally, for $y\in \D^*$ we have
\bel{TO5}
\Tilde\Omega(y)~=~\Tilde \tau(y)\,e^{\tau(y)} h(\xi(y))
 + e^{\tau(y)}\,\nabla h(\xi(y))\cdot \tilde\xi(y) .\eeq
  In view of (\ref{ttau})--(\ref{tx2}) and the definition of $\Tilde C$ at (\ref{TC}), this yields
 \bel{TO6}\bega{rl}\ds
 \bigl|\Tilde\Omega(y)\bigr|&\ds\leq~e^{\tau(y)}\|h\|_{\L^\infty}
\cdot{\bigl|\bfw(\tau(y),y)\bigr|\over\left| \la \bfn(\xi(y)))\,,~\bfq(\xi(y)) \ra\right|} \\[4mm]
&\ds\qquad\qquad   + 
e^{\tau(y)}\|\nabla h\|_{\L^\infty}\cdot   \left( 1+  { |\bfq(\xi(y))|\over \left| \la \bfn(\xi(y)),\,\bfq(\xi(y))
\ra\right| } \right)\,\bigl|\bfw(\tau(y),y)\bigr|\\[5mm]
&\leq~e^{\tau(y)}\Tilde C\,\bigl|\bfw(\tau(y),y)\bigr|~\leq~e^{T^*}\Tilde C\,\bigl|\bfw(\tau(y),y)\bigr|.\enda\eeq
\v
{\bf 2.}  
 It now remains to derive a bound on the $\L^2$ norm of 
$\bfw$.
Keeping (\ref{tT}) in mind, define
$$Z(t)~\doteq~\int_\D |\bfw(t, y)|^2\, dy.$$
For any $t\in [0, T^*]$,  by (\ref{qass})  
the Jacobian determinant of the map $y\mapsto x(t,y)$ satisfies
\bel{det}\det \Big({\partial x(t,y)\over\partial y}\Big) ~=~e^{2\mu t}.\eeq
Using the above identity to change variables of integration, 
setting
$$\D^*_t~\doteq~\bigl\{ y\in \D^*\,;~~\tau(y)<t\bigr\},$$
we obtain
\bel{i3}\int_{\D^*_t} |\Tilde\bfq (x(t,y))|^2\, dy~\leq~e^{-2\mu t} \int_{\D^*} |\Tilde\bfq (x)|^2\, dx~
=~e^{-2\mu t}\,\|\Tilde\bfq \|_{\L^2(\D^*)}^2\,.\eeq
In turn, by the elementary inequality $ab\leq {1\over 2} (a^2+b^2)$, this yields
$$\bega{rl}\ds {d\over dt} Z(t)&\ds\leq~2\int_{\D^*_t} \Big|\la \bfw(t,y), ~\dot\bfw(t,y)\ra\Big|\, dy
\\[4mm]
&\ds\leq~2\int_{\D^*_t}
 \Big\{ \bigl|D\bfq(x(t,y))\bigr|\, |\bfw(t,y)|^2 + |\bfw(t,y)|\, |\Tilde\bfq (x(t,y))|\Big\}\, dy
\\[4mm]
&\ds\leq~2 \Big\{ \|D\bfq\|_{\C^0} \cdot \|\bfw(t,\cdot)\|_{\L^2(\D^*_t)}^2  + \|\bfw(t,\cdot)\|_{\L^2(\D^*_t)}
e^{-\mu t}\|\Tilde \bfq\|_{\L^2(\D^*)}\Big\}\\[4mm]
&\leq~\ds 2 \Big\{ C\, Z(t)  + {Z(t)\over 2} + 
{e^{-2\mu t}\over 2} \|\Tilde\bfq \|^2_{\L^2(\D^*)}\Big\}
\,.
\enda$$
Hence
\bel{wl2}\|\bfw(t,\cdot)\|_{\L^2(\D)}^2~=~Z(t)~\leq~\kappa (t)\cdot \|\Tilde \bfq\|^2_{\L^2(\D)}, \eeq
where we set
\bel{kapdef} \kappa (t)~\doteq~{1\over 2}
 \int_0^t e^{(2M+1)(t-s)} \,e^{-2\mu s}\, ds~=~{e^{(2M+1)t}- e^{-2\mu t}\over 4M+4\mu+2}\,.\eeq
 Taking $t=T^*$ we conclude
\bel{wl3}\int_{\D^*} \bigl|\bfw(\tau(y),y)\bigr|^2\, dy~\leq~{e^{(2M+1)T^*}- e^{-2\mu T^*}\over 4M+4\mu+2}\cdot \|\Tilde \bfq\|^2_{\L^2(\D)}
\,.\eeq
Using this bound, from (\ref{TO6}) one obtains
(\ref{oper}). 
 \endproof

\subsection{A regularity estimate.}  
Given a H\"older continuous function
$f:\D\mapsto\R$, for $0<\alpha,\delta<1$, we introduce the notation
\bel{39}
\|f\|_{\alpha,\delta}~\doteq~
\sup_{x,y\in\D,~0<|x-y|<\delta} ~{|f(x)-f(y)|\over |x-y|^\alpha}\,.\eeq
Notice that, compared with the standard definition of the norm in the H\"older space $\C^{0,\alpha}$
(see for example \cite{BFA, Evans}), in (\ref{39}) the supremum is taken only over couples with $|x-y|<\delta$. 
From the above definition it immediately follows
\bel{ades}
\|f\|_{\alpha,\delta}~\leq~\|f\|_{\C^{0,\alpha}}\,,
\qquad\quad \|f\|_{\alpha,\delta}~\leq~
\delta^{1-\alpha} \|\nabla f\|_{\L^\infty}\,.
\eeq

The next lemma shows that the $\C^{0,\alpha}$ norm can be controlled in terms
of the seminorm 
$\|\cdot\|_{\alpha,\delta}$ together with the $\L^2$ norm.

As a preliminary we observe that, since $\D$ is an open set with  smooth boundary, it has a positive inner radius $\rho>0$.
Namely, every point $y\in \D$ lies in an open ball 
of radius $\rho$, entirely contained inside $\D$.

\begin{lemma}\label{l:interp}
Let $\D\subset\R^2$ be a bounded open set with inner radius $\rho>0$, and let $\alpha,\delta\in \,]0,1[\,$ and $c>0$ be given.   Then there exists $\delta_c>0$ such that the following holds.
If
\bel{40} \|f\|_{\alpha,\delta}~\leq ~c,\qquad\qquad \|f\|_{\L^2(\D)}~\leq~\delta_c\,,\eeq
then \bel{42}
\|f\|_{\C^0(\D)}~\leq~{c\over 2} \,\delta^\alpha.\eeq
In turn, this implies
\bel{41}\|f\|_{\C^{0,\alpha}(\D)}~\leq~c.
\eeq
\end{lemma}

{\bf Proof.} 
{\bf 1.} We claim that,  by choosing $\delta_c>0$ small enough, the inequalities
(\ref{40}) imply (\ref{42}).
Indeed, consider the function
$$\phi(x)~\doteq~\max\Big\{ 0,~{c\over 2} \delta^\alpha - c|x|^\alpha\Big\}.$$
and the disc
$$B~=~\{(x_1,x_2)\,;~~(x_1-\rho)^2 + x_2^2<\rho^2\}.$$
Setting 
\bel{doc}\delta_c~\doteq~\|\phi\|_{\L^2(B)} ~=~\left(\int_B\phi^2(x)\, dx\right)^{1/2},\eeq
we claim that the assumptions (\ref{40}) imply (\ref{42}).
Indeed, assume that  $f(y_0)> c\delta^\alpha/2$ at some point $y_0\in \D$.
Let $B_0\subset\D$ be a disc of radius $\rho$ which contains $y_0$.
A comparison argument now yields
$$\bega{rl}\ds\|f\|^2_{\L^2(\D)}&\ds\ge~\int_{B_0} f^2(x)\, dx~\geq~\int_{B_0}
\left(\max\Big\{ 0,~|f(y_0)|-c|x-y_0|^\alpha\Big\}\right)^2\, dx \\[4mm]
\ds\qquad &\ds >~
\int_B\phi^2(x)\, dx~=~\delta_c^2\,,\enda$$
reaching a contradiction.  Hence (\ref{42}) holds.
\v
{\bf 2.} By the assumptions, we already know that 
\bel{44} {|f(x)-f(y)|\over |x-y|^\alpha}~\leq~c\eeq
when $0<|x-y|<\delta$.   It remains to prove that the same holds when
$|x-y|\geq\delta$.   But in this case by (\ref{42}) one trivially has
$${|f(x)-f(y)|\over |x-y|^\alpha}~\leq~{|f(x)|+|f(y)|\over \delta^\alpha}
~\leq~c\,.$$
Together with (\ref{42}), this yields
$$\|f\|_{\C^{0,\alpha}(\D)}~\doteq~\max\left\{ \sup_x~|f(x)|\,,~~\sup_{x\not= y}  {|f(x)-f(y)|\over |x-y|^\alpha}\right\}
~\leq~c\,,$$
proving (\ref{41}).
\endproof

\section{Construction of an exact solution}
\setcounter{equation}{0}
\label{s:4}
As in (\ref{UV})-(\ref{proj}), we consider the decomposition $H\doteq \L^2(\D)=U\times V$, with perpendicular projections 
$P$ and $ I-P$, and write $\Omega= (u,v)$. We recall that $D_u\Lambda$, $D_v\Lambda$ denote the partial derivatives of the map
$\Lambda=\Lambda(u,v)$ introduced at (\ref{ldef}).

As anticipated in Section~\ref{s:2}, we write the fixed point problem $\Omega = \Lambda(\Omega)$ in the
equivalent form
\bel{ups}
\Omega~=~\Ups(\Omega)~\doteq~(I-AP)^{-1}\bigl(\Lambda(\Omega)- AP\,\Omega\bigr).\eeq
In terms of the components $(u,v)$, at 
$\Omega=\Omega_0= (u_0,0)$, the differential of the map $\Omega\mapsto \Lambda(\Omega)- AP\,\Omega$
has the form (\ref{diff}).    

To achieve the contraction property, on the product space $H=U\times V$  we consider 
the equivalent inner product 
\bel{ip*}\la(u,v),(u',v')\ra_*~=~\langle u, u'\rangle + \eta_0 \langle v, v'\rangle,\eeq
for a suitable constant $0<\eta_0\leq 1$.  The corresponding norm is
\bel{n*}\|(u,v)\|_*~=~\bigl( \|u\|^2 + \eta_0  \|v\|^2\bigr)^{1/2}.\eeq
Based on (\ref{dism}), we assume that a constant $\eta_0$ can be chosen so that, at the point 
$\Omega_0$, the corresponding norm of the linear operator (\ref{diff}) is $\leq 1/4\gamma$.
By (\ref{IA2}) and the definition of $\Ups$ at (\ref{ups}),  this implies
$$\|D\Ups(\Omega_0)\|_*~\leq~{1\over 4}\,.$$
Here and in the sequel, we also denote by $\|\cdot\|_*$ the norm of a linear operator, corresponding to the 
norm (\ref{n*}) on the product space $H=U\times V$.

By continuity, we can determine a radius $r_0>0$ such that,
denoting by $B_*(\Omega_0, r_0)$ a ball centered at $\Omega_0$ with radius $r_0$ 
w.r.t.~the equivalent norm $\|\cdot\|_*$, one has the implication
\bel{c23}
\Omega,\Omega'\in B_*(\Omega_0, r_0)\qquad \implies\qquad \|\Ups(\Omega)-\Ups(\Omega')\|_*~\leq~
{1\over 2} \|\Omega-\Omega'\|_*\,.\eeq
If the approximate solution $\Omega_0$ satisfies
\bel{ig}
\|\Ups(\Omega_0)-\Omega_0\|_*~\leq~{r_0\over 2}\,,\eeq
we can then define the iterates
\bel{OPN}
\Omega_n~\doteq~\Ups(\Omega_{n-1}),\qquad\qquad \Phi_n~=~\Delta^{-1}\Omega_n\,.\eeq
Taking the limit 
\bel{UPN}
\ov\Omega~=~\lim_{n\to\infty} \Omega_n\,,
\eeq
we thus obtain  a fixed point
$\ov\Omega = \Ups(\ov\Omega)$, with 
\bel{fixo}
\|\ov\Omega-\Omega_0\|_{\L^2(\D)}~\leq~\|\ov\Omega-\Omega_0\|_*~\leq~r_0\,.\eeq

While this approach is entirely straightforward, our main concern here is to 
derive more precise estimates on the various constants, which guarantee that an exact solution actually exists.

Notice that the contraction property, in a norm equivalent to $\L^2$, implies that 
$\|\Ups^n(\Omega_0)-\Omega_0\|_{\L^2}$ will be small, for every $n\geq 1$.   
In turn, recalling (\ref{OPN}), we conclude that
$\|\Phi_n-\Phi_0\|_{H^2}$  also remains small.    However, this estimate
is not enough to provide an a priori bound on $\|\nabla^\perp \Phi_n\|_{\C^1}$, which is needed to 
estimate the differential $D\Lambda$.  For this reason, an additional regularity estimate for the 
iterates $\Omega_n=\Ups^n(\Omega_0)$ will be derived.


\subsection{Regularity estimates.}
Assuming that the $\C^1$  vector field $\bfq_0= \nabla^\perp \Phi_0+\bfv$ satisfies the transversality 
assumptions {\bf (A1)}, we can find $0<\delta_1\leq 1$ with the following property.
\begi
\item[{\bf (P1)}]
{\it If $\|\Phi-\Phi_0\|_{\C^2}\leq\delta_1$, then the vector field $\bfq= \nabla^\perp \Phi+\bfv$
still satisfies (\ref{OT3})--(\ref{OT2}), possibly with a smaller constant $c_1>0$  and with $T^*$ replaced by $T^*+1$.

}
\endi
In particular, all characteristics starting from a point $y\in\Sigma_1$ in the support of $h$ still exit from the domain $\D$ 
within time $T^*+1$.
Moreover, the right hand side of (\ref{TC}) remains uniformly bounded by some constant, 
which we still denote by  $\Tilde C$.

Our next goal is to ensure that all our approximations  satisfy
\bel{pno}
\|\Phi_n-\Phi_0\|_{\C^2}~\leq~ \delta_1\qquad\qquad \forall n\geq 1.\eeq 
This bound will be achieved by an inductive argument, in several steps.
\v
{\bf 1.} 
If (\ref{pno}) holds, 
then 
\bel{qn}\|\bfq_n\|_{\C^1}~=~\bigl\|\nabla^\perp \Phi_n+ \bfv\bigr\|_{\C^1}~\leq~
\bigl\|\nabla^\perp \Phi_0+ \bfv\bigr\|_{\C^1} + \bigl\|\nabla^\perp \Phi_n- \nabla^\perp \Phi_0\bigr\|_{\C^1}
~\leq~M+1\,.\eeq
Assuming (\ref{pno}), the solution $\Omega = \Gamma(\Phi_n)$ of (\ref{lam1}) satisfies
\bel{Gn1}
\|\Omega\|_{\C^0}~\leq~e^{T^*+1}\, \|h\|_{\C^0}\,.\eeq
To provide a bound on the gradient $\nabla\Omega$, recalling  the notation introduced at 
(\ref{xty})--(\ref{tauy}), 
fix any $x_0\in \D$ such that $x_0 = \exp(t\bfq_n)(y)$ for some $y\in \Sigma_1\cap \Supp(h)$ and $t= \tau(x_0)>0$. 
Using the notation $x(t, x_0) \doteq \exp(-t\bfq_n)(x_0)$, consider the vector
$$\bfw(t)~\doteq~\lim_{\ve\to 0} {x(t, x_0+\ve \bfe) - x(t,x_0)\over\ve}\,,$$
where $\bfe\in \R^2$ is any unit vector.
Then 
$$\dot \bfw(t)~=~D\bfq_n(x(t, x_0))\cdot \bfw(t),\qquad\qquad \bfw(0)~=~\bfe.$$
Hence
\bel{wn}
|\bfw(t)|~\leq ~\exp\bigl\{ t\|D\bfq_n\|_{\C^0}\bigr\}~\leq~e^{ (T^*+1) (M+1)}.\eeq
By (\ref{OR}), the same computations performed at (\ref{wwy})--(\ref{TO6}) 
now yield
$$|\nabla \Omega(x_0)\cdot \bfe|~\leq~e^{T^*+1} \,\Tilde C\, \bigl|\bfw(\tau(x_0))\bigr|~\leq~\Tilde C\, e^{(T^*+1)(M+2)}.$$
Since the unit vector $\bfe$ was arbitrary, this implies that the gradient of $\Omega=\Gamma(\Phi_n)$ satisfies
\bel{nom}
\bigl\|\nabla \Omega\bigr\|_{\C^0}~\leq~\Tilde C\, e^{(T^*+1)(M+2)}.\eeq
\v
{\bf 2.}
Recalling (\ref{ups}) and the definition of $\Gamma$ at (\ref{lam1})-(\ref{ldef}), 
by (\ref{Gn1}) and (\ref{nom}) we now obtain
\bel{on1}
\bigl\|\Gamma(\Delta^{-1}\Omega_n)\bigr\|_{\C^1}~\leq~e^{T^*+1} \|h\|_{\C^0} +\Tilde C\, e^{(T^*+1)(M+2)},\eeq
as long as (\ref{pno}) holds.  In turn, this yields a bound of the form
\bel{oc1}
\|\Omega_{n+1}\|_{\C^1}~=~\Big\| (I-AP)^{-1} \bigl(\Gamma(\Delta^{-1}\Omega_n)- AP \,\Omega_n\bigr)\Big\|_{\C^1}\leq~C_1\,,\eeq
where the constant $C_1$ can be estimated in terms of (\ref{on1})
and the properties of the linear operator $A$.
\v
{\bf 3.} 
Since the domain $\D$ has smooth boundary, Schauder's regularity estimates \cite{Evans, GT} with $\alpha=1/2$
yield a bound of the form 
\bel{pnes}
\|\Phi_n-\Phi_0\|_{\C^2}~=~\|\Delta^{-1}(\Omega_n-\Omega_0)\|_{\C^{2}} ~\leq~C_2\, \|\Omega_n-\Omega_0\|_{\C^{0,
1/2}}\,,\eeq
for some constant $C_2$ depending only on $\D$.

Assume we have the inductive estimate
\bel{ie2}
\|\Omega_n-\Omega_0\|_{\C^1} ~\leq~2C_1\,.
\eeq
Choosing $0<\delta<1$ so that  $\delta\leq (2C_1 C_2)^{1\over\alpha-1} = (2 C_1C_2)^{-2}$,
we obtain
\bel{del}\|\Omega_n-\Omega_0\|_{\alpha,\delta} ~\leq~\delta^{1/2}\|\nabla\Omega_n-\nabla\Omega_0\|_{C^0}
~\leq~2C_1\delta^{1/2}~\leq~{1\over C_2} \,.\eeq
\v
{\bf 4.} Let $\delta_1>0$ be the constant introduced in {\bf (P1)}.
We now use Lemma~\ref{l:interp}, with $\delta>0$ as in (\ref{del}) and $c=\delta_1/C_2$.   
This yields  a constant $\delta_c\in \, ]0, \delta_1]$ such that
(\ref{40}) implies (\ref{42})-(\ref{41}).   In the present setting, this means that the two inequalities
\bel{ndc}\|\Omega_n-\Omega_0\|_{\C^1} ~\leq~2C_1\,,\qquad\qquad
 \|\Omega_n-\Omega_0\|_{\L^2} ~\leq~\delta_c\,,\eeq
together imply
\bel{onz}\|\Omega_n-\Omega_0\|_{\C^{0,1/2}}~\leq~2C_1\,\delta^{1/2} ~\leq~{\delta_1\over C_2}\,,
\qquad\qquad \|\Phi_n-\Phi_0\|_{\C^2}~\leq
~\delta_1.\eeq
In other words, if the $\L^2$ distance between $\Omega_0$ and every $\Omega_n$
remains small, then all these approximate solutions have uniformly bounded $\C^{0,1/2}$ norm.
In turn, property {\bf (P1)} applies.

\subsection{Convergence of the approximations.}

As long as $\|\Omega-\Omega_0\|_{\L^2}\leq \delta_c$ we have $\|\Phi\|_{\C^2}\leq \delta_1$, and hence
\bel{kod}\|D\Gamma(\Omega)\|_{\L^2}~\leq~\Tilde C\cdot K(M+\delta_1, \, T^*+1)~\doteq~\kappa_0\,.\eeq
\bel{kol}\|D\Lambda(\Omega)\|_{\L^2}~\leq~\|D\Gamma(\Omega)\|_{\L^2}
\cdot \|\Delta^{-1}\|_{\L^2} ~\leq~~\kappa_0/\lambda_1\,,
\eeq
where $\lambda_1>0$ denotes the first  eigenvalue of the Laplace operator on $\D$.

On the other hand, choosing a sufficiently large number $N$ of functions in the orthogonal basis 
at (\ref{UV}), we achieve (\ref{LVsmall}).  When $(u,v)=(u_0,0)\in U\times V$,
the four blocks in the matrix of partial derivatives (\ref{diff}) have norms which can be dominated
respectively by the entries of the matrix
$\left(\bega{cc} 0 & \kappa_0\ve_0\cr\kappa_0/\lambda_1 & \kappa_0\ve_0\enda\right)$.
More generally, we can determine $\delta_2\in \,]0, \delta_c]$ such that
\bel{imp3}
\|\Omega-\Omega_0\|_{\L^2}~\leq~\delta_2\qquad\implies\qquad \|P \,D_u \Lambda(\Omega)- P\, D_u\Lambda(\Omega_0)\|_{\L^2}
~\leq~\kappa_0\ve_0\,.\eeq
Computing the matrix of partial derivatives at (\ref{diff}) at such a point $\Omega=(u,v)$, and recalling
that $A\doteq P\, D_u\Lambda(\Omega_0)$, we obtain the relation 
\bel{cbo}
D\bigl(\Lambda(u,v)-Au\bigr)
~=~\left(\bega{ccc} P \,D_u \Lambda(\Omega)- P\, D_u\Lambda(\Omega_0) && P D_v\Lambda\\[4mm]
(I-P) D_u\Lambda && (I-P) D_v\Lambda\enda\right)~~\prec~~
\left(\bega{cc} \kappa_0\ve_0 & \kappa_0\ve_0\cr\kappa_0/\lambda_1 & \kappa_0\ve_0\enda\right).
\eeq
Here we have used the notation $A\prec B$, meaning that every entry in the $2\times 2$ matrix $A$ has norm
bounded by the corresponding entry in the matrix $B$.  Notice that the constant $\delta_2$ in (\ref{imp3})
 involves only
the behavior of $\Gamma\circ \Delta^{-1}$ on a neighborhood of $\Omega_0$ in the finite dimensional subspace
$U$, and can be directly estimated.

According to {\bf (A3)}, we now make the key assumption that the constant 
$\ve_0>0$  in (\ref{LVsmall}) is small enough so that
\bel{A3}\sqrt 2\, \kappa_0\left( \ve_0^2+{\ve_0\over\lambda_1}
\right)^{1/2}
\leq ~{1\over 2\gamma}\,.\eeq
This allows us to introduce on $U\times V$ the equivalent norm (\ref{n*}), choosing 
$$\eta_0~\doteq~\lambda_1\ve_0\,.$$
Without loss of generality, we can assume that $0<\eta_0\leq 1$, so that 
\bel{eqn}\|\cdot\|_*~\leq~\|\cdot\|_{\L^2}~\leq~{1\over\sqrt{\eta_0}} \|\cdot\|_*\,.\eeq
In term of this new norm, a direct computation shows, at any point $(u,v)$ where (\ref{cbo}) holds,
the corresponding operator norm satisfies\footnote{Indeed, assume $u,v\in \R$,  $u^2+\eta_0v^2\leq 1$.
Set $z=\sqrt{\eta_0} v$, so that $u^2+z^2\leq 1$.
This implies
$$\bega{l}\ds(\kappa_0 \ve_0 u  + \kappa_0\ve_0 v)^2 + \eta_0\left({\kappa_0\over \lambda_1} u +\kappa_0\ve_0 v
\right)^2~=~\left( \kappa_0 \ve_0 u  + {\kappa_0\ve_0 \over\sqrt {\lambda_1\ve_0 }}z\right)^2
+ \lambda_1\ve_0\left({\kappa_0\over \lambda_1} u +{\kappa_0\ve_0 \over\sqrt {\lambda_1\ve_0 }}z
\right)^2\\[4mm]
\ds\qquad =~\kappa_0^2\left[ \Big( \ve_0 u + \sqrt{\ve_0\over\lambda_1} z\Big)^2 + 
\Big(\sqrt{\ve_0\over\lambda_1}\, u
+ \ve_0 z\Big)^2\right]~=~\kappa_0^2\left[ \ve_0^2( u^2+z^2)  + {\ve_0\over\lambda_1}(u^2+z^2) 
 + 4\ve_0\sqrt{\ve_0\over\lambda_1}\, uz \right]\\[4mm]
 \qquad \leq~\ds 2\kappa_0^2 \left( \ve_0^2 + {\ve_0\over\lambda_1}\right) (u^2+z^2).
\enda$$
An entirely similar computation, with the numbers $u,v$ replaced by $\|u\|$ and $\|v\|$ respectively, shows that 
the corresponding norm of the linear operator (\ref{cbo}) is bounded by $\sqrt 2\, \kappa_0\left( \ve_0^2+\ds{\ve_0\over\lambda_1}
\right)^{1/2}$. }

\bel{opn3}\Big\|
D\bigl(\Lambda(u,v)-Au\bigr)\Big\|_*~\leq~{1\over 2\gamma}\,.\eeq
Hence
\bel{co2}\|D\Ups\|_*~\leq~\|(I-PA)^{-1}\|\cdot \Big\|D\bigl(\Lambda(u,v)-Au\bigr)\Big\|_*~\leq~{1\over 2}\,.\eeq

By the previous analysis, restricted to the set
$${\cal S}~\doteq~\left\{ \Omega\,;~~\|\Omega-\Omega_0\|_{\L^2}~\leq~ \delta_2\,,\quad \|\Omega-\Omega_0
\|_{\C^{0,\alpha}}\,\leq\, {\delta_1\over C_2}
\right\},$$
the map $\Omega\mapsto \Ups(\Omega)$ is a strict contraction, w.r.t.~the equivalent norm $\|\cdot \|_*$ introduced at 
(\ref{n*}).  Indeed, by (\ref{co2}),
\bel{co1}
\Omega,\Omega'\,\in\, {\cal S}\qquad\implies\qquad \|\Ups(\Omega)-\Ups(\Omega')\|_*
~\leq~{1\over 2} \|\Omega-\Omega'\|_*\,.\eeq

Recalling that $\delta_1$ is the constant in {\bf (P1)}, we now have the chain of implications
$$\bega{l}\ds \|\Omega-\Omega_0\|_*~\leq~ {\delta_2\over\sqrt{\eta_0}} \quad\implies\quad
 \|\Omega-\Omega_0\|_{\L^2}
~\leq~\delta_2\\[4mm]\ds
\quad\implies\quad \|\Omega-\Omega_0
\|_{\C^{0,\alpha}}\,\leq\, {\delta_1\over C_2}\quad\implies\quad \|\Delta^{-1}(\Omega-\Omega_0)\|_{\C^2}
~\leq~\delta_1\,.\enda$$

If  the initial guess $\Omega_0$ satisfies 
\bel{igss}\|\Ups(\Omega_0) - \Omega_0\|_{\L^2}~\leq~{\delta_2\sqrt{\eta_0}\over 2}\,,\eeq
then by (\ref{eqn}) and (\ref{co1}) all the iterates $\Omega_n$ in (\ref{OPN}) will remain inside ${\cal S}$.
Namely,
$$\|\Omega_n-\Omega_0\|_{\L^2}~\leq~{1\over\sqrt{\eta_0}}\,\|\Omega_n-\Omega_0\|_*~\leq~
{1\over\sqrt{\eta_0}}\cdot  2\|\Ups(\Omega_0)-\Omega_0\|_*~\leq~\delta_2\,.$$
Letting $n\to \infty$ we have the convergence $\Omega_n\to \ov\Omega\in {\cal S}$.
This limit function $\ov\Omega$ provides a solution to the boundary value problem 
(\ref{lam1}), (\ref{D-1}), with 
\bel{cl1}\|\ov\Omega-\Omega_0\|_{\L^2}~\leq~\delta_2\,.\eeq

\subsection{An existence theorem.}   For readers' convenience, 
we summarize the previous analysis, recalling the various constants introduced along the way.
All these constants can be estimated in terms of the domain $\D$, the boundary data $g,h$ in (\ref{BC}),
and the finite dimensional
approximation $\Omega_0$ produced by the numerical algorithm. 
\begi
\item $\Tilde C$ is the constant at (\ref{TC}), related to the stability of the ODEs for the characteristics of the linear PDE
(\ref{lam3}).  This applies to any vector field $\bfq= \nabla^\perp\Phi+\bfv$ with $\|\Phi-\Phi_0\|_{\C^2}\leq \delta_1$.

\item $K(M,t)$ is the function in (\ref{pera}).   Together with $\Tilde C$, it provides an estimate (\ref{oper}) on 
how the solution  to the linear PDE (\ref{lam3}) varies in the $\L^2$ distance, 
depending on the vector field $\nabla^\perp \Phi$.

\item $C_2$ is the Schauder regularity constant for the smooth domain $\D$, introduced at (\ref{pnes}).

\item $\lambda_1>0$ is the lowest eigenvalue of the Laplace operator on $\D$.

\item $\ve_0$ is the small constant in (\ref{LVsmall}), determining the rate at which the components in the orthogonal 
space $V$ are damped by the inverse Laplace operator $\Delta^{-1}$.

\item $\kappa_0$ is the constant in (\ref{kod}), estimating how the solution of the linear PDE (\ref{lam1})
varies, depending on the vector field $\nabla^\perp\Phi_1$, w.r.t.~the $\L^2$ norm.

\item $\delta_1$ is the constant in {\bf (P1)}, measuring by how much we can perturb the vector field
$\bfq= \nabla^\perp\Phi+\bfv$, and still retain the transversality conditions at the boundary, and a finite exit time.

\item $\delta_c$ is the constant in the (\ref{ndc}), providing the regularity estimate (\ref{onz}).

\item $\gamma$ is a bound on the norm $\|(I-PD_u\Lambda(\Omega_0))^{-1}\|$ of the inverse Jacobian matrix in (\ref{IA2}).
It gives a measure of stability for the finite dimensional fixed point problem $\Omega = P\circ
\Lambda(\Omega)$.

\item $\delta_2$ is the constant in (\ref{imp3}), determining a neighborhood of the initial approximation $\Omega_0$,
where the differential $D\Lambda(\Omega)$ remains close to $D\Lambda(\Omega_0)$. 
\endi

Recalling the definition of the finite dimensional linear operator $A$ at  (\ref{Adef}) and the equivalent formulation
of the boundary value problem (\ref{lam2})-(\ref{lam1}) as a fixed point problem (\ref{ups}), 
we can now summarize all of the previous analysis as follows. 

\begin{theorem}  Let $\D\subset\R^2$ be a bounded open domain with smooth boundary, decomposed as 
$\partial\D=\Sigma_1\cup\Sigma_2$.   Given smooth boundary data $g,h$,  consider the boundary value problem 
(\ref{lam2})-(\ref{lam1}), where $\bfv$ is the vector field in (\ref{bbv}).  
Assume that the properties {\bf (A1)-(A2)} hold, together with {\bf (P1)}.

Consider the orthogonal decomposition $\L^2(\D)= U\times V$ as in (\ref{UV})-(\ref{proj}), 
and let
 $\Omega_0= (u_0, 0)\in U\times V$  be an approximate solution such that  
\bel{igs2}\|\Ups(\Omega_0) - \Omega_0\|_{\L^2(\D)}~\leq~{\delta_2\sqrt{\lambda_1\ve_0}\over 2}\,.\eeq
Then an exact solution $\ov\Omega$ exists, with
\bel{OO} 
\|\ov\Omega-\Omega_0\|_{\L^2(\D)}~\leq~\delta_2\,.
\eeq

\end{theorem}

\v
{\bf Acknowledgment.} The authors would like to thank Ludmil Zikatanov for useful discussions.
\v

\end{document}